\documentclass[11pt,a4paper]{article}

\usepackage[latin1]{inputenc}
\usepackage{amssymb}
\usepackage{amsmath}
\usepackage{stmaryrd}
\usepackage[english]{babel}
\usepackage{cases}
\usepackage{fancyhdr}
\usepackage{graphicx}
\usepackage{color}
\usepackage[latin1]{inputenc}
\usepackage{color}
\usepackage{latexsym}
\usepackage{amssymb}
\usepackage{graphicx}
\usepackage{enumerate}
\usepackage{extarrows}
\newtheorem{Theorem}{Theorem}[part]

\newtheorem{Proposition}{Proposition}[part]

\newtheorem{Lemma}{Lemma}[part]

\newenvironment{proof}[1][Proof]{\textbf{#1.} }{\ \rule{0.5em}{0.5em}}

\def \2{\vspace{2mm}}

\def \nn {\nonumber}

\def \Sum{\displaystyle\sum}

\def \Frac{\displaystyle\frac}

\def \R{\mathbb{R}}

\def \E{\mathbb{E}}

\def \Ec{{\cal E}}

\def \Hf    {\mathfrak{H}}

\def \o {\otimes }

\def \ep{\hbox{ }\hfill$\Box$}

\def\be{\begin{eqnarray}}
\def\ee{\end{eqnarray}}

\def\b*{\begin{eqnarray*}}
\def\e*{\end{eqnarray*}}
\def \nn{\nonumber }

\begin{document}

\begin{center}
{\large{\bf Berry-Esséen bounds and almost sure CLT for the quadratic variation of the bifractional Brownian motion}}\\~\\
Soufiane Aazizi$^*$\footnote{Department of Mathematics, Faculty of
Sciences Semlalia Cadi Ayyad University, B.P. 2390 Marrakesh,
Morocco. Email: {\tt aazizi.soufiane@gmail.com}

$^*$This author is supported by the Marie Curie Initial Training
Network (ITN) project: ``Deterministic and Stochastic Controlled
Systems and Application", FP7-PEOPLE-2007-1-1-ITN, No. 213841-2.}
and Khalifa Es-Sebaiy\footnote{ENSA de Marrakech, Université Cadi
Ayyad, Marrakech, Maroc. Email: {\tt k.Essebaiy@uca.ma} }
\\
{\it  Université Cadi Ayyad}\\~\\
\end{center}

%

\begin{abstract} ~~\\Let $B$ be a bifractional Brownian motion with parameters
$H\in (0, 1)$ and $K\in(0,1]$.  For any $n\geq1$, set $Z_n =
\sum_{i=0 }^{n-1}\big[
n^{2HK}(B_{(i+1)/n}-B_{i/n})^2-\E\left((B_{i+1}-B_{i})^2\right)\big]$.
We use the Malliavin calculus and the so-called Stein's method on
Wiener chaos introduced by Nourdin and Peccati \cite{NP09} to
derive, in the case when $0<HK\leq3/4$,
  Berry-Ess\'een-type bounds for the Kolmogorov distance between the
law of the correct renormalization $V_n$ of $Z_n$   and the standard
normal law.  Finally, we study almost sure central limit theorems
for the sequence  $V_n$.
\end{abstract}

\vspace{13mm}

\noindent {\bf Key words~:} Kolmogorov distance; Central limit
theorem; Almost sure central limit theorem;
 Bifractional Brownian motion; Multiple stochastic integrals; Quadratic variation.

\vspace{5mm}


\vspace{13mm}

\section{Introduction}
\setcounter{equation}{0} \setcounter{Assumption}{0}
\setcounter{Theorem}{0} \setcounter{Proposition}{0}
\setcounter{Corollary}{0} \setcounter{Lemma}{0}
\setcounter{Definition}{0} \setcounter{Remark}{0}

Let $B = ( B_t, t\geq 0 )$ be a  bifractional Brownian motion
(bifBm)  with parameters $H\in (0, 1)$ and $K\in(0,1]$, defined on
some probability space $(\Omega,\mathcal{F},P)$. (Here, and
everywhere else, we do assume that $\mathcal{F}$ is the sigma-field
generated by $B$.) This means that $B$ is a centered
Gaussian process with the covariance function $E[B_sB_t]=R_{H,K}(s,t)$, where%
\begin{align}
\label{cov}
R_{H,K}(s,t)=\frac{1}{2^K}\Big(\big(t^{2H}+s^{2H}\big)^K-|t-s|^{2HK}\Big).
\end{align}
  The case $K = 1$ corresponds to the fractional Brownian
motion (fBm) with Hurst parameter H. The process $B$ has no
stationary increments,   but it has the quasi-helix property (in the
sense of J.P. Kahane), \be \label{quasi-helix}
2^{-K}|t-s|^{2HK}\leq\E\left(\left|B_t-B_s\right|^2\right)\leq
2^{1-K}|t-s|^{2HK},\ee so $B$ has $\gamma-$H\"{o}lder continuous
paths for any $\gamma\in (0,HK)$ thanks to the Kolmogorov-Centsov
theorem, and it is a self-similar process, that is, for any constant
$a>0$, the processes $(B_{at},\ t\geq0)$ and $(a^{HK}B_{t},\
t\geq0)$ have the same distribution. The bifBm $B$ can be extended
for $1 < K < 2$ with $H\in (0, 1)$ and $HK\in (0, 1)$ (see
\cite{Bar-Ess11}). We refer to \cite{HV03, RT06, ET07, LN09} for
further details on the subject.

 An example of
interesting problem related to $B$ is the study of the asymptotic
behavior of the quadratic variation of $B$ on $[0,1]$ defined as
\[Z_n =
\sum_{i=0 }^{n-1}\left[
n^{2HK}(B_{(i+1)/n}-B_{i/n})^2-\E\left((B_{i+1}-B_{i})^2\right)\right],
\quad n\geq1.\] Let us consider the correct renormalization $V_n$ of
$Z_n$ given as, \be \label{V_n}V_n=\frac{Z_n}{\sqrt{Var(Z_n)}}.\ee
Recall that, if $Y$, $Z$ are two real-valued random variables, then
the Kolmogorov distance between the law of  $Y$ and the law of $Z$
is given by
\begin{eqnarray*}d_{\mbox{\tiny{Kol}}}(Y,Z)=\sup_{-\infty<z<\infty}|P(Y\leq z)-P(Z\leq z)|.
\end{eqnarray*}
In the particular case of the fBm  (that is when $K =1$), and thanks
to  the seminal works of Breuer and Major \cite{BM83}, Dobrushin and
Major \cite{DM79}, Giraitis and Surgailis \cite{GS85} and Taqqu
\cite{T79}, it is well-known that we have, as $n\rightarrow\infty$:
\begin{itemize}
\item If $0<H<\frac34$ then
\[\frac{V_n}{\sigma_H\sqrt{n}}\overset{\rm
law}{\longrightarrow} \mathcal{N }(0, 1).\]
\item If $H=\frac34$ then\[\frac{V_n}{\sigma_H\sqrt{n\log(n)}}\overset{\rm
law}{\longrightarrow} \mathcal{N }(0, 1).\]
\item If $H>\frac34$ then\[\frac{V_n}{n^{2H-1}}\overset{\rm
law}{\longrightarrow} Z\sim \mbox{``Hermite random variable" }.\]
\end{itemize}Here, $\sigma_H> 0$ denotes an (explicit) constant depending only on $H$. Moreover, explicit bounds for the Kolmogorov distance between the law of
 $V_n$    and the standard normal
law are obtained by \cite[Theorem 4.1]{NP09}, \cite[Theorem
1.2]{BN08} and \cite[Theorem 5.6]{Nourdin-Prix}. The following facts
happen: For some constant $c_H$ depending only on $H$, we have:
 \b*
    d_{Kol}\left(V_n,\mathcal{N}(0, 1)\right) \le c_{H}\times
    \left\{
     \begin{array}{ll}
       \frac{1}{\sqrt{n}}   & \mbox{ if } H\in \left(0, \frac58\right)  \\
       ~~\\ \frac{(\log n)^{3/2}}{\sqrt{n}}   &\mbox{ if } H=\frac58  \\
       ~~\\
       n^{4H-3}  &\mbox{ if } H  \in \left(\frac58 ,\frac{3}{4}\right) \\
       ~~\\
       \frac{1}{\sqrt{\log n}}  &\mbox{ if } H  =\frac{3}{4}
 \end{array}
   \right.
\e* On other hand, Bercu et al. \cite{BNT10}  proved the almost sure
central limit theorem  (ASCLT) for $V_n$. Recently, Tudor \cite{T11}
studied the subfractional Brownian motion case.

Let us now describe the results we prove in the present paper.
First, in Theorem \ref{ThemBerryEsseen} we use   the Malliavin
calculus and Stein method, in the case when $HK\in(0,\frac{3}{4}]$,
to derive explicit bounds for the Kolmogorov distance between the
law of $V_n$    and the standard normal law. Precisely, three cases
are considered according to the value of $HK$:
 \b*
    d_{Kol}\left(V_n,\mathcal{N}(0, 1)\right) \le c_{H,K}\times
    \left\{
     \begin{array}{ll}
       n^{-\frac12}   &\mbox{ if } HK\in \left(0, \frac12\right] \\
       ~~\\
       n^{2HK-\frac{3}{2}} &\mbox{ if } HK  \in \left[\frac12 ,\frac{3}{4}\right)\\
       ~~\\
       \frac{1}{\sqrt{\log n}} &\mbox{ if } HK  =\frac{3}{4}
 \end{array}
   \right.
\e*  where $c_{H,K}$ is a constant depending only on $H$ and $K$. In
 Theorem \ref{ASCLT V_n}, we prove
almost sure central limit theorem for  $V_n$.

The rest of the paper is organized as follows. Section 2  deals with
preliminaries concerning Malliavin calculus, Stein's method and
related topics needed throughout the paper. Section 3 and 4 contain
our main results, concerning Berry-Essén bounds and ASCLT for the
quadratic variation of the bifractional Brownian motion.

\section{Preliminaries}
In this section, we briefly recall some basic facts concerning
Gaussian analysis and Malliavin calculus that are used in this
paper; we refer to \cite{N06} for further details. Let $\Hf$ be a
real separable Hilbert space. For any $q \geq 1$, we denote by
$\Hf^{\otimes q}$ (resp. $\Hf^{\odot q}$) the $q$th tensor product
(resp. $q$th symmetric tensor product) of $\Hf$. We write $X =
\{X(h), h\in \Hf\}$ to indicate a centered isonormal Gaussian
process on $\Hf$. This means that $X$ is a centered Gaussian family,
defined on some probability space $(\Omega ,\mathcal{F}, P)$ and
such that $E [X(g)X(h)] = \langle g, h\rangle_{\Hf}$ for every $g,
h\in\Hf$. (Here, and everywhere else, we do assume that
$\mathcal{F}$ is the sigma-field generated by $X$.)

For every $q\geq 1$, let $\mathcal{H}_q$ be the $q$th Wiener chaos
of $X$, that is, the closed linear subspace of $L^{2}(\Omega)$
generated by the random variables $\{H_{q}\left( X\left( h\right)
\right) ,h\in \Hf,\| h\| _{\Hf}=1\}$, where $H_{q}$ is the $q$th
Hermite polynomial defined as $H_{q}(x)=(-1)^q
e^{\frac{x^2}{2}}\frac{d^q}{dx^q}(e^{-\frac{x^2}{2}})$. The mapping
$I_{q}(h^{\otimes q})=H_{q}\left( X\left( h\right) \right) $
provides a linear isometry between the symmetric tensor product
$\Hf^{\odot q}$ (equipped with the modified norm
$\|\cdot\|_{\Hf^{\odot q}}= \sqrt{q!}  \|\cdot\|_{\Hf^{\otimes q}}$)
and $\mathcal{H}_q$. Specifically, for all $f,g\in\Hf^{\odot q}$ and
$q\geq 1$, one has
\begin{equation}\label{isoint}
E\big[I_q(f)I_q(g)\big]=q!\langle f,g\rangle_{\mathcal{H}^{\otimes
q}}
\end{equation}
On the other hand, it is well-known  that any random variable $Z$
belonging to $L^2(\Omega)$ admits the following chaotic expansion:
\begin{equation}\label{chaos}
Z=E[Z]+\sum_{q=1}^\infty I_q(f_q)
\end{equation}
where the series converges in $L^2(\Omega)$ and the kernels $f_q$,
belonging to $\Hf^{\odot q}$,
 are
uniquely determined by $Z$.

Let $\{e_k, k\geq 1\}$ be a complete orthonormal system in $\Hf$.
Given $f\in\Hf^{\odot p}$ and $g\in\Hf^{\odot q}$, for every $r=0,
\dots, p\wedge q$, the $r$th contraction of $f$ and $g$ is the
element of $\Hf^{\otimes(p+q-2r)}$ defined as
$$
f\otimes_r g=\sum_{i_1=1, \dots, i_r=1}^{\infty} \langle
f,e_{i_1}\otimes\cdots\otimes e_{i_r}\rangle_{\Hf^{\otimes r}}
\otimes \langle g,e_{i_1}\otimes\cdots\otimes
e_{i_r}\rangle_{\Hf^{\otimes r}}.
$$
In particular, note that $f\otimes_0g=f\otimes g$ and when $p=q$,
that $f\otimes_pg=\langle f, g\rangle_{\Hf^{\otimes p}}$.
 Since, in
general, the contraction $f\otimes_rg$ is not necessarily symmetric,
we denote its symmetrization by $f\widetilde\otimes_rg \in
\Hf^{\odot(p+q-2r)}$. When $f\in \Hf^{\odot q}$, we write $I_q(f)$
to indicate its $q$th multiple integral with respect to $X$. The
following formula is useful to compute the product of such multiple
integrals: if $f\in \Hf^{\odot p}$ and $g\in \Hf^{\odot q}$, then
\begin{equation}
\label{eq:multiplication} I_p(f)I_q(g)=\sum_{r=0}^{p\wedge q} r!
\left(\!\!\begin{array}{c}p\\r\end{array}\!\!\right)
\left(\!\!\begin{array}{c}q\\r\end{array}\!\!\right)
I_{p+q-2r}(f\widetilde\otimes_rg).
\end{equation}

 Let $\cal S$ be the set of all smooth cylindrical random
variables, that is, which can be expressed as $F = f(X(\phi_1),
\ldots, X(\phi_n))$ where $n\geq 1$, $f : \R^n \rightarrow \R$ is a
$\mathcal{C}^\infty$-function such that $f$ and all its derivatives
have at most polynomial growth, and $\phi_i\in\Hf$. The Malliavin
derivative of $F$ with respect to $X$ is the square integrable
$\Hf$-valued random variable defined as
\[
DF\; =\; \sum_{i =1}^n \frac{\partial f}{\partial x_i}(X(\phi_1),
\ldots, X(\phi_n)) \phi_i.
\]
In particular, $DX(h) = h$ for every $h\in \Hf$. As usual,
$\mathbb{D}^{1,2}$ denotes the closure of the set of smooth random
variables with respect to the norm
$$\| F\|_{1,2}^2 \; = \; E[F^2] +
E\big[\|DF\|_{\Hf}^2\big].$$ The Malliavin derivative $D$ verifies
the chain rule: if $\varphi:\R^n\rightarrow\R$ is $\mathcal{C}^1_b$
and if $(F_i)_{i=1,\ldots,n}$ is a sequence of elements of
$\mathbb{D}^{1,2}$, then $\varphi(F_1,\ldots,G_n)\in
\mathbb{D}^{1,2}$ and we have
$$
D\varphi(F_1,\ldots,G_n)=\sum_{i=1}^n
\frac{\partial\varphi}{\partial x_i} (F_1,\ldots, G_n)DF_i.
$$

Recall the following results concerning CLT and ASCLT for multiple
stochastic integrals.
\begin{Theorem}[Nourdin-Peccati \cite{NP09}]\label{Kol}  Let $q \ge 2$ be an integer and let $F=I_q(f)$ with  $f\in{{\cal{H}}^{\odot
q}}$, then
\begin{eqnarray}\label{Kol expression}d_{\mbox{\tiny{Kol}}}(F,N)
\leq\sqrt{E\left[\left(1-\frac{1}{q}\|DF\|_{{\cal{H}}}^2\right)^2\right]},
\end{eqnarray}where $N\sim\mathcal{N}(0,1)$.
\end{Theorem}
\begin{Theorem}[Bercu et al. \cite{BNT10}]\label{ThemBercu}Let $q \ge 2$ be an integer, and let $\{G_n\}_{n\geq1}$  be a
sequence of the form $G_n=I_q(f_n)$, with $f_n\in{{\cal{H}}^{\odot
q}}$. Assume that $E[G^2_n] = q!\|f_n\|^2_{\Hf^{\o q}}  = 1$ for all
$n$, and  that $G_n\overset{\rm law}{\longrightarrow} N\sim
\mathcal{N }(0, 1)$ as $n\rightarrow\infty$. If the two following
conditions are satisfied
\begin{enumerate}[1)]
  \item $
        \Sum_{n=2}^\infty \frac{1}{n \log^2n}\Sum_{k=1}^n\frac{1}{k}\|f_k\o_r f_k\|_{\Hf^{\o 2(q-r)}} <
        \infty$ for avery $1 \le r \le q-1$,

  \item
       $ \Sum_{n=2}^\infty \frac{1}{n \log^3n}\Sum_{k,l=1}^n\frac{\left|\langle f_k,f_l\rangle_{\Hf^{\o q}}\right|}{kl}<
       \infty,$
\end{enumerate}
then $\{G_n\}_{n\geq1}$ satisfies an ASCLT. In other words, almost
surely, for any bounded and continuous function
$\varphi:\R\rightarrow\R,$
\[\frac{1}{\log(n)}\sum_{k=1}^{n}\frac{1}{k}\varphi(G_k)\rightarrow\E\varphi(N),\quad \mbox{as }n\rightarrow\infty.\]
\end{Theorem}

From now, assume on one hand that $X=B$ is a bifBm  with parameters
$H\in(0,1)$ and $K\in(0,1]$ and on the other hand $\Hf $ is a real
separable Hilbert space defined as follows: (i) denote by $\Ec$ the
set of all $\R$-valued step functions on $[0,\infty)$, (ii) define
$\Hf$ as the Hilbert space obtained by closing $\Ec$ with respect to
the scalar product
 $$\langle1_{[0, s]}, 1_{[0, t]}\rangle_\Hf=R_{H,K}(s,t)=\frac{1}{2^K}\Big(\big(t^{2H}+s^{2H}\big)^K-|t-s|^{2HK}\Big).$$
 In particular,  one has that $B_t=B(1_{[0,t]} )$.
\section{Berry-Esséen bounds in the CLT for the quadratic variation of the bifBm}
\setcounter{equation}{0} \setcounter{Assumption}{0}
\setcounter{Theorem}{0} \setcounter{Proposition}{0}
\setcounter{Corollary}{0} \setcounter{Lemma}{0}
\setcounter{Definition}{0} \setcounter{Remark}{0}

In this section, we  prove that for every  $HK \in \left(0,
\frac{3}{4}\right]$ a Central Limit Theorem holds, where $V_n$ was
defined in (\ref{V_n}). Using the Stein's method we also derive the
Berry-Esséen bounds for this convergence.
\subsection{General setup}
Let us define
 \b*
\theta(i,j)=2^{-K}\big(\gamma(i,j)+\rho(i-j)\big),\quad i,j\in
\mathbb{N} \e* where \be
\gamma(i,j)&=&\left((i+1)^{2H} +(j+1)^{2H}\right)^K -\left(i^{2H} +(j+1)^{2H}\right)^K  - \left((i+1)^{2H} + j^{2H}\right)^K \nn \\
&&+\left(i^{2H} + j^{2H}\right)^K, \label{gamma}\ee and \be
\rho(r)&=&|r+1|^{2HK}+|r-1|^{2HK}-2|r|^{2HK},\quad r\in \mathbb{Z
}.\ee
~~\\
Observe that the function $\gamma$ is symmetric, $\rho(0)=2$,
$\rho(x)=\rho(-x)$ and $\rho$ behaves asymptotically as
\be\label{behavior rho} \rho(r)=2HK(2HK-1)|r|^{2HK-2},\quad
|r|\rightarrow\infty.\ee In particular,
$\sum_{r\in\mathbb{Z}}\rho^2(r)<\infty$ if, and only if,
$HK\in(0,\frac34)$.\\
 We will use the notation \be\delta_{k/n}=1_{[k/n,(k+1)/n]}&\mbox{ and
 }&
 \sigma= \sqrt{\frac{1}{8}\sum_{r\in \mathbb{Z}}\rho^2(r)}.\label{segma and delta_k/n}\ee
Using self-similarity property of $B$ and (\ref{cov}) we deduce that
 \b*n^{2HK}\langle\delta_{i/n},\delta_{j/n}\rangle_\Hf&=&
  n^{2HK}\E\left(\left(B_{\frac{i+1}{n}}-B_{\frac{i}{n}}\right)\left(B_{\frac{j+1}{n}}-B_{\frac{j}{n}}\right)\right)\\
 &=&\E\left( (B_{i+1}-B_i)(B_{j+1}-B_j)\right)
 \\&=&\theta(i,j).
 \e*
Hence, we can write the quadratic variation of $B$, with respect to
a subdivision $\mathcal{\pi}_n = \{0< \frac{1}{n}
<\frac{2}{n}<\ldots < 1\}$ of $[0, 1]$, as follows
 \be
Z_n&=&\Sum_{k=0}^{n-1}\left[ n^{2HK}\left(B_{\frac{k+1}{n}}-B_{\frac{k}{n}}\right)^2 -\theta(k,k)\right]\nn\\
&=&\Sum_{k=0}^{n-1}\left[ n^{2HK}\left(I_1(\delta_{k/n})\right)^2 -\theta(k,k)\right]\nn\\
&=&I_2\left(\underbrace{n^{2HK} \Sum_{k=0}^{n-1}\delta^{\o 2}_{k/n}}_{g_n}\right)\nn\\
&=&I_2(g_n). \label{Z_n} \ee Thus, we can also write the correct
renormalization $V_n$, defined in (\ref{V_n}), of $Z_n$  as follows,
\be \label{V_n second
expression}V_n=\frac{Z_n}{\sqrt{Var(Z_n)}}=\frac{I_2(g_n)}{\sqrt{Var(Z_n)}}.\ee
 Before computing the Kolmogorov
distance, we start with the following results which are used
throughout the paper. Here, and everywhere else, the notation $a_n
\trianglelefteqslant b_n$ means that $\sup_{n\ge1} |a_n|/|b_n| <
\infty$.

\begin{Lemma}\label{Lemmagamma}
~~
\begin{enumerate}[i)]
  \item  Fixing $y\geq0$ (resp. $x\geq 0$), the function $x\rightarrow\gamma(x,y)$ (resp. $y\rightarrow\gamma(x,y)$) defined in (\ref{gamma}) is
  increasing
   for $H \in \left(0, \frac{1}{2}\right]$.
  \item  For any $H \in (0, 1)$ and $K\in (0, 1]$, the function $\gamma$ is negative and we
   have for $j$ large
    \be
    \gamma(0,j) &\sim& c_{H,K} j^{2HK-2}\label{gamma(0,j)},\\
    \gamma(j,j)&\sim& c_{H,K} j^{2HK-2}\label{gamma(j,j)}.
    \ee

    If $j\leq l$ then
    \be
|\gamma(j,l)| \le c_{H,K} l^{2HK-2}.   \label{|gamma(k,l)|}
    \ee
     where $c_{H,K}$ is a constant (explicit) depending only on $H$ and $K$.
\end{enumerate}

\end{Lemma}
\proof  $i)$  We fix $y\geq0$, \be
\frac{\partial \gamma}{\partial x}(x,y)&=&2HK(x+1)^{2H-1}\Big[((x+1)^{2H} +(y+1)^{2H})^{K-1} - ((x+1)^{2H} + y^{2H})^{K-1}\Big] \nn\\
      &&-2HKx^{2H-1}\Big[(x^{2H} +(y+1)^{2H})^{K-1} -(x^{2H} + y^{2H})^{K-1} \Big]\nn\\
      &=&2HK\Big[g(1+x)-g(x) \Big],
\ee where \b* g(x)=x^{2H-1}\Big[(x^{2H} +(y+1)^{2H})^{K-1} -(x^{2H}
+ y^{2H})^{K-1}\Big], \e* If $H \in \left(0, \frac{1}{2}\right]$ and
$K\in (0, 1]$, then $\gamma$ is increasing since the function $g$ is
increasing on $(0, \infty)$. Indeed, \b*
g'(x)&=&(2H-1)x^{2H-2}\Big[(x^{2H} +(y+1)^{2H})^{K-1} -(x^{2H} + y^{2H})^{K-1}\Big]\\
    &&+2H(K-1)x^{4H-2}\Big[(x^{2H} +(y+1)^{2H})^{K-2} -(x^{2H} + y^{2H})^{K-2}\Big]\\
    &\ge& 0.
\e*
 $ii)$ To show that $\gamma$ is negative, it suffices to
remark the decreasing property of the function
  $p:x\in[0,\infty) \rightarrow (a+x)^K - (b+x)^K$.\\
   By straightforward expansion of function $\gamma$, we can easily
prove $(\ref{gamma(0,j)})$ and $(\ref{gamma(j,j)})$.
\\
   If $H \leq
\frac{1}{2}$, by the first point $i)$, the function
$x\rightarrow|\gamma(x,y)|$  is
  decreasing. Thus, we deduce
\b* |\gamma(k,l)| \le |\gamma(0,l)| \sim c_{H,K} l^{2HK-2}. \e* If
$H
> \frac{1}{2}$, we rewrite $\gamma$ as $\gamma(k,l)
=g_k(1+l)-g_k(l)$ where $g_k(x):=  ((k+1)^{2H} + x^{2H})^K-(k^{2H} +
x^{2H})^K$. Applying  mean value theorem we obtain for some $x_{k,l}
\in [l, l+1]$ that\b* |\gamma(k,l)|&=&2HKx_{k,l}
^{2H-1}\Big[\big(x_{k,l} ^{2H} +k^{2H}\big)^{K-1} -\big(x_{k,l}
^{2H} + (k+1)^{2H}\big)^{K-1}\Big]\\&\leq&2HK(l+1)
^{2H-1}\Big[\big(l ^{2H} +k^{2H}\big)^{K-1} -\big(l ^{2H} +
(k+1)^{2H}\big)^{K-1}\Big]. \e* Again by mean value theorem on
$y\rightarrow\big(l^{2H} + y^{2H}\big)^{K-1}$, we have for some
$y_{k,l} \in [k, k+1]$ \b* \Big[\big(l ^{2H} +k^{2H}\big)^{K-1}
-\big(l ^{2H} +
(k+1)^{2H}\big)^{K-1}\Big]=2H(K-1)y_{k,l}^{2H-1}\Big[l^{2H}
+y_{k,l}^{2H}\Big]^{K-2}. \e* Consequently, for $k\leq l$, \b*
|\gamma(k,l)|&\leq&4H^2K(1-K)(l+1)^{2H-1}(k+1)^{2H-1}\Big[l^{2H}
+k^{2H}\Big]^{K-2}\\&\leq &c_{H,K}l^{2HK-2} \e*  and the second
point $ii)$ follows.\ep
\begin{Proposition} \label{ProplimitVn} Let $Z_n$ be the sequence  defined in
(\ref{Z_n}) and  let $\sigma$ be the constant given by (\ref{segma
and delta_k/n}).
\begin{itemize}\item[1.] Assume that $0 < HK< \frac{3}{4}$. Then, as $n\rightarrow\infty$, it holds
\be \frac{Var(Z_n)}{4^{2-K}n\sigma^2} \longrightarrow 1
\label{limitVn}. \ee
\item[2.] Assume that $HK= \frac{3}{4}$. Then, as $n\rightarrow\infty$, it holds
\be \frac{Var(Z_n)}{4^{2-K}\sigma^2n\log n} \longrightarrow 1
\label{limitVnLogn}. \ee
\end{itemize}
\end{Proposition}
\proof To show (\ref{limitVn}), we write \b*
\frac{Var(Z_n)}{4^{2-K}n\sigma^2}-1&=&\Frac{n^{-1}}{4^{2-K}\sigma^2}~\E[I_2^2(g_n)]-1=\frac{n^{-1}}{2^{3-2K}\sigma^2}~\|g_n\|^2_{\Hf^{\o 2}}-1\\
    &=&\frac{n^{4HK-1}}{2^{3-2K}\sigma^2}  \Sum_{k,l=0}^{n-1} \langle \delta^{\o 2}_{k/n},\delta^{\o 2}_{l/n}\rangle_{\Hf^{\o 2}}-1\\
    &=&\frac{n^{4HK-1}}{2^{3-2K}\sigma^2}   \Sum_{k,l=0}^{n-1} \langle \delta_{k/n},\delta_{l/n}\rangle_{\Hf}^2-1
    \e*
    \b*
    &=&\frac{n^{-1}}{2^{3-2K}\sigma^2}  \Sum_{k,l=0}^{n-1} \theta^2(k,l)-1\\
    &=&\frac{n^{-1}}{8\sigma^2}  \Sum_{k,l=0}^{n-1} \gamma^2(k,l)
    +\left(\frac{n^{-1}}{8\sigma^2}  \Sum_{k,l=0}^{n-1} \rho^2(k-l)-1\right) +\frac{n^{-1}}{4\sigma^2}  \Sum_{k,l=0}^{n-1}\gamma(k,l)\rho(k-l)\\
    &=:&J_1(n)+J_2(n)+J_3(n).
\e* As in the proof of  \cite[Theorem 4.1]{NP09}, we have \be
\label{J_2(n)}J_2(n)\longrightarrow0,\quad\mbox{ as }
n\rightarrow\infty.\ee
 On the other hand
 \b*J_1(n)&=&\frac{n^{-1}}{8\sigma^2}  \Sum_{k,l=0}^{n-1} \gamma^2(k,l)
 \\&=&\frac{n^{-1}}{8\sigma^2}  \Sum_{k=0}^{n-1}
 \gamma^2(k,k)+\frac{n^{-1}}{4\sigma^2}  \Sum_{0\leq k<l\leq n-1} \gamma^2(k,l)\\
 &=:&J_{1,1}(n)+J_{1,2}(n).
 \e*
By (\ref{gamma(j,j)}), the sum
\b*J_{1,1}(n)=\frac{n^{-1}}{8\sigma^2} \Sum_{k=0}^{n-1}
 \gamma^2(k,k),
 \e*
 behaves as $\Frac{n^{-1}}{8\sigma^2} \;\sum_{k=0}^{n-1} k^{4HK-4}$ which goes to zero as $n\rightarrow\infty$, because
 $HK<3/4$.\\Thus,
\be \label{J_{1,1}} J_{1,1}(n)\longrightarrow0\quad\mbox{as
}n\rightarrow\infty.
 \ee
Now, we study the convergence of $J_{1,2}(n)$. We first fix two
positive constants $\alpha$ and $\beta$ such that $\alpha+\beta=1$
and $4HK-2<\beta<1$. \vspace{2mm}
\\
We deduce from (\ref{|gamma(k,l)|}) that \b*
J_{1,2}(n)=\frac{n^{-1}}{4\sigma^2} \Sum_{0\leq k<l\leq n-1}
\gamma^2(k,l) &\le&C_{H,K}\frac{n^{-1}}{4\sigma^2}\Sum_{0\le l\le
n-1}l^{4HK-3}\\&\leq&
C_{H,K}\Frac{n^{-\alpha}}{4\sigma^2}~\Sum_{0\le l\le
n-1}l^{4HK-3-\beta}\longrightarrow0, \mbox{ as }n\rightarrow\infty.
\e* Hence, \be \label{J_{1,2}}
J_{1,2}(n)\longrightarrow0,\quad\mbox{as }n\rightarrow\infty.
 \ee Combining (\ref{J_{1,1}}) and (\ref{J_{1,2}}) leads to
 \be \label{J_1(n)}
J_{1}(n)\longrightarrow0,\quad\mbox{as }n\rightarrow\infty.
 \ee
Finally, from (\ref{J_1(n)}) and (\ref{J_2(n)})  together with
Cauchy Schwartz inequality \be
|J_3(n)|&\leq&\frac{n^{-1}}{4\sigma^2}  \Sum_{k,l=0}^{n-1}|\gamma(k,l)\rho(k-l)|\nn\\
   &\le&   \left(\frac{n^{-1}}{4\sigma^2}\Sum_{k,l=0}^{n-1}\gamma^2(k,l)\right)^{1/2}
\left(  \frac{n^{-1}}{4\sigma^2}\Sum_{k,l=0}^{n-1}\rho^2(k-l)\right)^{1/2}\nn\\
&=&2\sqrt{J_1(n)(J_2(n)+1)} \longrightarrow 0, \mbox{ as
}n\rightarrow\infty,\label{J_3(n)} \ee and the convergence
(\ref{limitVn}) follows.\\
We prove now (\ref{limitVnLogn}).  Following similar argument of the
proof of (\ref{limitVn}), we have \b*
\frac{Var(Z_n)}{4^{2-K}\sigma^2n\log n}-1
    &=&\frac{n^{-1}}{8\sigma^2\log n}  \Sum_{k,l=0}^{n-1} \gamma^2(k,l) +\left(\frac{n^{-1}}{8\sigma^2\log n}
    \Sum_{k,l=0}^{n-1} \rho^2(k-l)-1\right)\\
    &&+\frac{n^{-1}}{4\sigma^2\log n}  \Sum_{k,l=0}^{n-1}\gamma(k,l)\rho(k-l)\\
    &=&\frac1{\log n}J_1(n)+\frac1{\log n}J_2(n)+\frac1{\log n}J_3(n).
\e* From \cite[page 490]{BN08} we have \be
\label{J'_2(n)}\frac1{\log n}J_2(n)\longrightarrow0\quad\mbox{ as }
n\rightarrow\infty.\ee
 On the other hand, since $HK=\frac34$ and the fact that $\log(n)\sim \sum_{1}^{n-1}\frac1k$ we
 deduce easily  from (\ref{J_1(n)}) and (\ref{J_3(n)})  that
 \b*
 \frac1{\log n}J_1(n) + \frac1{\log n}J_3(n) \underset{n \rightarrow \infty}{\longrightarrow 0}.
 \e*
 \ep \vspace{2mm}
\\
\subsection{A Berry-Esséen bound for  $0< HK\leq \frac34$}

 Our first main result is summarized in the following Theorem.
\begin{Theorem}\label{ThemBerryEsseen}
Let $N \sim \mathcal{N}(0, 1)$ and let $V_n$ be defined by (\ref{V_n
second expression}). Then  $V_n$ converges in distribution to $N$.
In addition, for some constant $c_{H,K}$ depending uniquely on $H$
and $K$, we have: for every $n \ge 1$,
  \b*
    d_{Kol}\left(V_n,N\right) \le c_{H,K}\times
    \left\{
     \begin{array}{ll}
      \frac{1}{\sqrt{ n}}  &\mbox{ if } HK\in \left(0, \frac12\right] \\
       ~~\\
       n^{2HK-\frac{3}{2}} &\mbox{ if } HK  \in \left[\frac12 ,\frac{3}{4}\right)\\
       ~~\\
       \frac{1}{\sqrt{\log n}} &\mbox{ if } HK  =\frac{3}{4}
 \end{array}
   \right.
\e*

\end{Theorem}
\proof From (\ref{Z_n}), we have \b*
DZ_n&=&2n^{2HK}\Sum_{k=0}^{n-1}I_1(\delta_{k/n})\delta_{k/n}, \e*
then \b*
\|DZ_n\|^2_\Hf&=&4n^{4HK}\Sum_{k,l=0}^{n-1}I_1(\delta_{k/n})I_1(\delta_{l/n})\langle
\delta_{k/n},\delta_{l/n}\rangle_\Hf, \e* by the multiplication
formula (\ref{eq:multiplication}), we get \b*
\|DZ_n\|^2_\Hf&=&4n^{4HK}\Sum_{k,l=0}^{n-1}I_2(\delta_{k/n}\widetilde{\o}
\delta_{l/n})\langle \delta_{k/n},\delta_{l/n}\rangle_\Hf
+ 4n^{4HK}\Sum_{k,l=0}^{n-1}\langle \delta_{k/n},\delta_{l/n}\rangle^2_\Hf\\
&=&4n^{4HK}\Sum_{k,l=0}^{n-1}I_2(\delta_{k/n}\widetilde{\o}
\delta_{l/n})\langle \delta_{k/n},\delta_{l/n}\rangle_\Hf+ \E
\|DZ_n\|^2_\Hf . \e* Combining this with the fact that $\E
\|DZ_n\|^2_\Hf=2Var(Z_n)$, we obtain \b* \frac{1}{2}\|DV_n\|^2_\Hf
-1
&=&\frac{2n^{4HK}}{Var(Z_n)}\Sum_{k,l=0}^{n-1}I_2(\delta_{k/n}\widetilde{\o}
\delta_{l/n})\langle \delta_{k/n},\delta_{l/n}\rangle_\Hf. \e*
 It follows that
\be &\E&\left[\left(\frac{1}{2}\left\|DV_n\right\|^2_\Hf -
1\right)^2\right]
\nn\\
&=& \frac{4n^{8HK}}{Var^2(Z_n)}
\E\left[\left( ~~\Sum_{k,l=0}^{n-1}I_2(\delta_{k/n}\widetilde{\o }\delta_{l/n})\langle \delta_{k/n},\delta_{l/n}\rangle_\Hf\right)^2\right]\nn\\
&=&\frac{8n^{8HK}}{Var^2(Z_n)}\Sum_{i,j,k,l=0}^{n-1}\langle
\delta_{i/n},\delta_{j/n}\rangle_\Hf \langle
\delta_{k/n},\delta_{l/n}\rangle_\Hf  \langle
\delta_{i/n}\widetilde{\o}\delta_{j/n},\delta_{k/n}\widetilde{\o}\delta_{l/n}\rangle_{\Hf^{\o
2}}. \label{n^2A_n} \\
&=&\frac{8n^2}{Var^2(Z_n)}A(n)\label{DV_n and A(n)}\ee
  where \b*
A(n)&=&n^{8HK-2}\Sum_{i,j,k,l=0}^{n-1}\langle
\delta_{i/n},\delta_{j/n}\rangle_\Hf ~\langle
 \delta_{k/n},\delta_{l/n}\rangle_\Hf ~\langle \delta_{i/n}\widetilde{\o}\delta_{j/n},\delta_{k/n}\widetilde{\o}\delta_{l/n}\rangle_{\Hf^{\o 2}}\nn\\
&=&\frac{n^{8HK-2}}{2}\Sum_{i,j,k,l=0}^{n-1}\langle \delta_{i/n},\delta_{j/n}\rangle_\Hf ~ \langle \delta_{k/n},\delta_{l/n}\rangle_\Hf ~\Big(\langle \delta_{i/n} ,\delta_{k/n}\rangle_{\Hf} ~ \langle \delta_{j/n} ,\delta_{l/n}\rangle_{\Hf}\nn\\
&&\hspace{65mm}+\langle \delta_{i/n} ,\delta_{l/n}\rangle_{\Hf}\langle \delta_{j/n} ,\delta_{k/n}\rangle_{\Hf}\Big)\nn\\
 &=&n^{8HK-2}\Sum_{i,j,k,l=0}^{n-1}\langle \delta_{i/n},\delta_{j/n}\rangle_\Hf
  ~ \langle \delta_{i/n},\delta_{k/n}\rangle_\Hf ~\langle \delta_{k/n} ,\delta_{l/n}\rangle_{\Hf} ~ \langle \delta_{j/n} ,\delta_{l/n}\rangle_{\Hf}
  \e*Hence, using that fact that for every $a, b\in\R$; $|ab|\leq\frac{1}{2}(a^2+b^2)$, we
  have
  \be|A(n)|
 &\leq&\frac{n^{8HK-2}}{2}\Sum_{i,j,k=0}^{n-1}\left|\langle \delta_{i/n},\delta_{j/n}\rangle_\Hf
  ~ \langle \delta_{i/n},\delta_{k/n}\rangle_\Hf\right| \left(\Sum_{l=0}^{n-1}\langle \delta_{k/n} ,\delta_{l/n}\rangle_{\Hf}^2
  \right)\nn\\&&+\frac{n^{8HK-2}}{2}\Sum_{i,j,k=0}^{n-1}\left|\langle \delta_{i/n},\delta_{j/n}\rangle_\Hf
  ~ \langle \delta_{i/n},\delta_{k/n}\rangle_\Hf\right| \left(\Sum_{l=0}^{n-1}\langle \delta_{j/n} ,\delta_{l/n}\rangle_{\Hf}^2
  \right)\nn\\
&=&n^{8HK-2}\Sum_{i,j,k=0}^{n-1}\left|\langle
\delta_{i/n},\delta_{j/n}\rangle_\Hf
  ~ \langle \delta_{i/n},\delta_{k/n}\rangle_\Hf\right| \left(\Sum_{l=0}^{n-1}\langle \delta_{k/n} ,\delta_{l/n}\rangle_{\Hf}^2
  \right)\label{|A(n)|}
  \ee
By (\ref{|gamma(k,l)|}) and (\ref{behavior rho}), we obtain
 \be n^{4HK}\Sum_{l=0}^{n-1}\langle
 \delta_{k/n},\delta_{l/n}\rangle_{\Hf}^2&=&\Sum_{l=0}^{n-1}\theta^2(k,l)\nn\\
&\leq&2^{1-2K}\left(\Sum_{l=0}^{n-1}\gamma^2(k,l)+\Sum_{l=0}^{n-1}\rho^2(k-l)\right)\nn\\
&=&2^{1-2K}\left(\Sum_{l=0}^{k}\gamma^2(k,l)+\Sum_{l=k+1}^{n-1}\gamma^2(k,l)+\Sum_{r=-k}^{n-1-k}\rho^2(r)\right)
\nn\\
&\leq&2^{1-2K}\left(\Sum_{l=0}^{k}k^{4HK-4}+\Sum_{l=1}^{n-1}l^{4HK-4}+2\Sum_{r=0}^{n-1}\rho^2(r)\right)
\nn\\
&\trianglelefteqslant&1+\Sum_{l=0}^{n-1}l^{4HK-4}.\label{majoration
sum delta_k delta_l} \ee

On the other hand, by using (\ref{behavior rho}) \be
 n^{4HK-2}\Sum_{i,j,k=0}^{n-1}\left|\langle
\delta_{i/n},\delta_{j/n}\rangle_\Hf
  ~ \langle \delta_{i/n},\delta_{k/n}\rangle_\Hf\right|&=&\frac{1}{n^2}\Sum_{i,j,k=0}^{n-1}\left|\theta(i,j)\theta(i,k)\right|
  \nn\ee \be
  &=&\frac{1}{n^2}\Sum_{i=0}^{n-1}\left(\Sum_{j=0}^{n-1}\left|\theta(i,j)\right|\right)^2
  \nn\\
  &\leq&\frac{1}{n^2}\Sum_{i=0}^{n-1}
  \left(\Sum_{j=0}^{n-1}\left|\gamma(i,j)\right|+\Sum_{j=0}^{n-1}\left|\rho(i-j)\right|\right)^2
 \nn\\&=&2^{-2K}\frac{1}{n^2}\Sum_{i=0}^{n-1}
  \left(\Sum_{j=0}^{i}\left|\gamma(i,j)\right|+\Sum_{j=i+1}^{n-1}\left|\gamma(i,j)\right|+\Sum_{r=-i}^{n-1-i}\left|\rho(r)\right|\right)^2
\nn\\
&\leq&2^{-2K}\frac{1}{n^2}\Sum_{i=1}^{n-1}
  \left(i^{2HK-1}+\Sum_{j=1}^{n-1}j^{2HK-2}+2\Sum_{r=0}^{n-1}\left|\rho(r)\right|\right)^2
\nn \\&\trianglelefteqslant&
 \frac{1}{n^2}\Sum_{i=1}^{n-1}i^{4HK-2}
  +\frac{1}{n}\left(\Sum_{j=1}^{n-1}j^{2HK-2}\right)^2.\label{D_0(n)}\ee
  By (\ref{|A(n)|}), (\ref{majoration
sum delta_k delta_l}) and (\ref{D_0(n)}),
  \be|A(n)|&\trianglelefteqslant&\frac{1}{n^2}\Sum_{i=1}^{n-1}i^{4HK-2}
  +\frac{1}{n}\left(\Sum_{j=1}^{n-1}j^{2HK-2}\right)^2\nn\\
  &:=&D(n).\label{D(n)}
  \ee
If $0<HK<\frac12$, \be D(n)&=&\frac{1}{n^2}\Sum_{i=1}^{n-1}i^{4HK-2}
  +\frac{1}{n}\left(\Sum_{j=1}^{n-1}j^{2HK-2}\right)^2
  \nn\\&\leq &\frac{1}{n}\Sum_{i=1}^{\infty}i^{4HK-3}
  +\frac{1}{n}\left(\Sum_{j=1}^{\infty}j^{2HK-2}\right)^2
  \nn\\&\trianglelefteqslant&\frac{1}{n}.\label{D(n) for HK<1/2}\ee
If $\frac12\leq HK<\frac34$, then, by using the fact that for all
$\alpha>-1$; $\sum_{k=1}^{n-1}r^{\alpha}\sim
n^{\alpha+1}/(\alpha+1)$ as $n\rightarrow\infty$, \be
D(n)&=&\frac{1}{n^2}\Sum_{i=1}^{n-1}i^{4HK-2}
  +\frac{1}{n}\left(\Sum_{j=1}^{n-1}j^{2HK-2}\right)^2
  \nn\\&\leq &\Sum_{i=1}^{n-1}i^{4HK-4}
  +\left(\Sum_{j=1}^{n-1}j^{2HK-\frac52}\right)^2
  \nn\\&\trianglelefteqslant&n^{4HK-3}.\label{D(n) for 3/4>HK>1/2}\ee
Combining (\ref{Kol expression}), (\ref{DV_n and A(n)}),
(\ref{limitVn}), (\ref{D(n) for HK<1/2}) and (\ref{D(n) for
3/4>HK>1/2}), we deduce that for every $0<HK<\frac34$,
 \b*
    d_{Kol}\left(V_n,N\right) \trianglelefteqslant
    \left\{
     \begin{array}{ll}
       \frac{1}{\sqrt{ n}}  &\mbox{ if } HK\in \left(0, \frac12\right] \\
       ~~\\
       n^{2HK-\frac{3}{2}} &\mbox{ if } HK  \in \left[\frac12
       ,\frac{3}{4}\right)
 \end{array}
   \right.
\e* Assume now that $HK=\frac34$. From (\ref{|A(n)|}),
(\ref{majoration sum delta_k delta_l}) and (\ref{D_0(n)})  together
with the fact that $\sum_{r=1}^{n-1}r^{-1}\sim \log(n)$ as
$n\rightarrow\infty$,
  \be \frac{|A(n)|}{\log^2(n)}&\trianglelefteqslant&\frac{1}{\log(n)}\left(\frac{1}{n^2}\Sum_{i=1}^{n-1}i^{-1}
  +\frac{1}{n}\left(\Sum_{j=1}^{n-1}j^{-\frac12}\right)^2\right)\nn\\
  &\trianglelefteqslant&\frac{1}{\log(n)} \label{A(n)/log2(n)},
  \ee
and this completes the proof of Theorem \ref{ThemBerryEsseen}. \ep
\section{Almost sure central limit Theorem}
\setcounter{equation}{0} \setcounter{Assumption}{0}
\setcounter{Theorem}{0} \setcounter{Proposition}{0}
\setcounter{Corollary}{0} \setcounter{Lemma}{0}
\setcounter{Definition}{0} \setcounter{Remark}{0} We are going now
to prove the second main result of this paper, which state the ASCLT
of the bifractional Brownian motion and its quadratic variation.
\begin{Proposition}\label{ASCLT for bifBm}For all $H\in(0,1)$ and $K\in(0,1]$, we have, almost surely, for any bounded and continuous function
$\varphi:\R\rightarrow\R,$
\[\frac{1}{\log(n)}\sum_{k=1}^{n}\frac{1}{k}\varphi(k^{-HK}B_k)\rightarrow\E\varphi(N),\quad \mbox{as }n\rightarrow\infty,\]
where $N\sim\mathcal{N}(0,1)$.
\end{Proposition}
\proof The proof is straightforward by applying \cite[Thorem 4.1 and
Corollary 3.7]{BNT10} and the fact that  \b*
|E[B_jB_l]|&=&2^{-K}\left((j^{2H}+l^{2H})^K-|j-l|^{2HK}\right)\\
&\le & 2^{-K}\left(j^{2HK}+l^{2HK}-|j-l|^{2HK}\right)\\
&=&2^{1-K}|E[B^{HK}_jB^{HK}_l]|, \e* where $B^{HK}$ is a fractional
Brownian motion with Hurst parameter $HK$.  \ep
\begin{Theorem}\label{ASCLT V_n}If $HK \in \left( 0, \frac{3}{4} \right]$, then the sequence $( V_n)_{n\ge 0}$ satisfies the ASCLT.
In other words, almost surely, for any bounded and continuous
function $\varphi:\R\rightarrow\R,$
\[\frac{1}{\log(n)}\sum_{k=1}^{n}\frac{1}{k}\varphi(V_k)\rightarrow\E\varphi(N),\quad \mbox{as }n\rightarrow\infty,\]
where $N\sim\mathcal{N}(0,1)$.
\end{Theorem}
\proof We shall make use of  Theorem \ref{ThemBercu}. From Theorem
\ref{ThemBerryEsseen},  $(V_n)_n$ satisfies the CLT, so that, it
remain to check conditions  $1)$ and $2)$. The cases $HK \in
\left(0, \frac{3}{4} \right)$ and $H=\frac34$ are treated
separately. By (\ref{V_n second expression}), we can write
$V_n=I_2(g_n)$ where
\begin{flushleft}
$g_n=\frac{n^{2HK}}{\sqrt{Var(Z_n)}}\Sum_{k=1}^n\delta^{\o 2}_{k/n},
$
\end{flushleft}
which implies that\b* g_n \otimes_1 g_n=
\frac{n^{4HK}}{Var(Z_n)}\Sum_{k,l=1}^{n}\langle\delta_{k/n},\delta_{l/n}\rangle_\Hf
~\delta_{k/n} \o \delta_{l/n}. \e* We deduce that \be \|g_n
\otimes_1 g_n\|^2_{\Hf^{\o 2}}=
\frac{n^2}{Var^2(Z_n)}A(n).\label{fnotimesfn} \ee Assume  that $HK
\in \left(0, \frac{3}{4} \right)$. Combining (\ref{limitVn}),
(\ref{D(n)}), (\ref{D(n) for HK<1/2}) and (\ref{D(n) for
3/4>HK>1/2}), we have \b*\|g_n \otimes_1 g_n\|^2_{\Hf^{\o
2}}\trianglelefteqslant( n^{-1}+n^{4HK-3}) \trianglelefteqslant
    \left\{
     \begin{array}{ll}
       n^{-1}  &\mbox{ if } HK\in \left(0, \frac{1}{2}\right) \\
        ~~\\
       n^{4HK-3} &\mbox{ if } HK  \in \left[\frac{1}{2}  ,\frac{3}{4}
       \right)
     \end{array}
   \right.
\e*
Consequently, condition $1$) in Theorem \ref{ThemBercu}  is satisfied.  \\
On the other hand, by (\ref{limitVn}), we have for $k <l$ \b*
\langle g_k,g_l\rangle_{\Hf^{\o
2}}&=&\frac{(kl)^{2HK}}{\sqrt{Var(Z_k)}\sqrt{Var(Z_l)}}\Sum_{i=0}^{k-1}
\Sum_{j=0}^{l-1}\langle \delta_{i/k},\delta_{j/l}\rangle^2_\Hf \\
&\le&c_{H,K}\frac{1}{\sqrt{kl}}\Sum_{i=0}^{k-1}\Sum_{j=0}^{l-1}\theta^2(i,j)\\
&\le&
c_{H,K}\frac{1}{\sqrt{kl}}\left[\Sum_{i=0}^{k-1}\Sum_{j=0}^{l-1}
\rho^2(i-j)+\left(\Sum_{0 \le i \leq j \le k-1}
+\Sum_{i=0}^{k-1}\Sum_{j=k}^{l-1} \right)\gamma^2(i,j)\right]
 \e*
As in the proof of \cite[Theorem 5.1, page 1621]{BNT10}, we obtain
that \b* \frac{1}{\sqrt{kl}}\Sum_{i=0}^{k-1}\Sum_{j=0}^{l-1}
\rho^2(i-j) \le c_{H,K} \sqrt{\frac{k}{l}}. \e* Using Lemma
\ref{Lemmagamma},
 we obtain
 \b*
\frac{1}{\sqrt{kl}}\Sum_{0\le i\le j \le k-1} \gamma^2(i,j) &\le&
c_{H,K}\sqrt{\frac{k}{l}}\Sum_{0\le i \le k-1}i^{4HK-4} \le
c_{H,K}}\sqrt{\frac{k}{l}. \e*
 Again from
Lemma \ref{Lemmagamma}, we have \b*
\frac{1}{\sqrt{kl}}\Sum_{i=0}^{k-1}\Sum_{j=k}^{l} \gamma^2(i,j)
&\le& \frac{1}{\sqrt{kl}}\Sum_{i=0}^{k-1}\Sum_{j=1}^{l} j^{4HK-4}\le
c_{H,K}\sqrt{\frac{k}{l}}. \e* Combining all the above bounds we
obtain \b* \langle f_k,f_l\rangle_{\Hf^{\o 2}} \le
c_{H,K}\sqrt{\frac{k}{l}}. \e*  Finally, condition $2)$ in Theorem
\ref{ThemBercu} is satisfied. \vspace{2mm}
\\
Now, suppose that $HK=\frac{3}{4}$. It follows from
(\ref{fnotimesfn}), (\ref{limitVnLogn}) and (\ref{A(n)/log2(n)}) \b*
\|g_k \otimes_1 g_k\|^2_{\Hf^{\o 2}}=
\frac{k^2\log^2k}{Var^2(Z_k)}\frac{A(k)}{\log^2(k)} \le c_{H,K}
\log^{-1} k. \e* Leads to \b*
        \Sum_{n=2}^\infty \frac{1}{n \log^2n}\Sum_{k=1}^n\frac{1}{k}\|g_k\o g_k\|_{\Hf^{\o 2}}
&\le&   c_{H,K} \Sum_{n=2}^\infty \frac{1}{n \log^2n}\Sum_{k=1}^n\frac{1}{k \sqrt{\log k}}\\
&\le&   c_{H,K}\Sum_{n=2}^\infty \frac{1}{n \log^{\frac{3}{2}}n}<
\infty. \e* To close the proof, it suffices  to show that \be
\langle g_k,g_l\rangle_{\Hf^{\o 2}}\le c_{H,K}\sqrt{\frac{k\log l}{l
\log k}}, && \forall k >l. \label{kloglllogk} \ee According to
(\ref{limitVnLogn}), we have \b* \langle g_k,g_l\rangle_{\Hf^{\o
2}}&=&\frac{(kl)^{2HK}}{\sqrt{Var(Z_k)}\sqrt{Var(Z_l)}}\Sum_{i=0}^{k-1}\Sum_{j=0}^{l-1}
\langle \delta_{i/k},\delta_{j/l}\rangle^2_\Hf \\
&\le&\frac{c_{H,K}}{\sqrt{l\log k}\sqrt{ k\log l}}\Sum_{i=0}^{k-1}\Sum_{j=0}^{l-1}\theta^2(i,j)\\
&\le& \frac{c_{H,K}}{\sqrt{l\log k}\sqrt{ k\log
l}}\left[\Sum_{i=0}^{k-1}\Sum_{j=0}^{l-1} \rho^2(i-j)+\left(\Sum_{0
\le i \leq j \le k-1} +\Sum_{i=0}^{k-1}\Sum_{j=k}^{l-1}
\right)\gamma^2(i,j)\right]
 \e*
As in the proof of \cite[Proposition 6.4, page 1625]{BNT10}, we have
 for
all $1 \le k \le l$
 \b*
 \frac{1}{\sqrt{l\log k}\sqrt{ k\log l}}\Sum_{i=0}^{k-1}\Sum_{j=0}^{l-1} \rho^2(i-j) \le c_{H,K}\sqrt{\frac{k\log l}{l \log k}}.
 \e*
Using Lemma \ref{Lemmagamma} and the fact that $\sum_{r=1}^{n-1}
r^{-1} \sim \log n$ as $n\rightarrow\infty$, we deduce that \be
\frac{1}{\sqrt{l\log k }\sqrt{k\log l }}\Sum_{0 \le i \le  j \le
k-1}\gamma^2(i,j) &\le& c_{H,K} \frac{k\log l}{\sqrt{l\log k
}\sqrt{k\log l }}\le c_{H,K}\sqrt{\frac{k\log l}{l \log k}}. \ee
Again from Lemma \ref{Lemmagamma}, we obtain \be
\frac{1}{\sqrt{l\log k }\sqrt{k\log l
}}\Sum_{i=0}^{k-1}\Sum_{j=k}^{l-1}\gamma^2(i,j) &\le& c_{H,K}
\sqrt{\frac{k\log l}{l \log k}} \ee  which complete the proof of the
Theorem \ref{ASCLT V_n}. \ep

\end{document}